\documentclass[letterpaper,11pt,twoside,reqno,nonamelimits]{amsart}
\usepackage[pagebackref,colorlinks=true,urlcolor=blue,linkcolor=blue,citecolor=blue]{hyperref}
\usepackage[top=1.65in,bottom=1.7in,left=1.5in,right=1.5in]{geometry}
\usepackage{times}
\usepackage{amsmath}
\usepackage{amssymb}
\usepackage{amsfonts}
\usepackage{amscd}
\usepackage{amsthm}
\usepackage{amsxtra}

\usepackage[all]{xy}  \CompileMatrices
\usepackage{mathrsfs}
\usepackage{nicefrac}
\usepackage{bbold}
\usepackage{graphicx}
\usepackage{color}
\usepackage{enumerate}
\usepackage{url}

\usepackage{algorithm}

\theoremstyle{plain}
\newtheorem{theorem}{Theorem}
\newtheorem{proposition}[theorem]{Proposition}
\newtheorem*{proposition*}{Proposition}

\newtheorem*{corollary*}{Corollary}
\newtheorem{lemma}[theorem]{Lemma}

\newtheorem*{theorem*}{Theorem}
\newtheorem*{lemma*}{Lemma}
\newtheorem*{conjecture*}{Conjecture}

\newtheorem*{question*}{Question}

\theoremstyle{definition}

\newtheorem*{exercise*}{Exercise}

\theoremstyle{remark}
\newtheorem{remark}[theorem]{Remark}
\newtheorem*{remark*}{Remark}
\newtheorem{remsTh}[theorem]{Remarks}

\newcommand{\subclass}[1]{}

\newcommand{\enumTi}[1]{\renewcommand{\theenumi}{#1}}

\newcommand{\alphenumi}{\enumTi{\alph{enumi}}}

\newcommand{\romenumi}{\enumTi{\roman{enumi}}}

\DeclareMathOperator{\rk}{rk}

\newcommand{\lt}{\left}
\newcommand{\rt}{\right}

\newcommand{\Tp}{{\mspace{-1mu}\scriptscriptstyle\top\mspace{-1mu}}}

\newcommand{\abs}[1]{{\lt\lvert{#1}\rt\rvert}}

\newcommand{\sabs}[1]{{\lvert{#1}\rvert}}

\newcommand{\widebar}[1]{\overline{#1}}

\newcommand{\nfrac}[2]{{\nicefrac{#1}{#2}}}

\newcommand{\RR}{\mathbb{R}}

\newcommand{\kk}{\mathbb{k}}

\newcommand{\Nm}[1]{{\lt\lVert #1 \rt\rVert}}

\newcommand{\ip}[1]{\lt(#1\rt)}

\newlength{\algotabbingwidth}
\DeclareMathOperator{\bc}{bc}
\DeclareMathOperator{\fool}{fool}
\DeclareMathOperator{\jk}{jk}

\begin{document}
\title[Biclique covering]{On some lower bounds on the number of bicliques needed to cover a bipartite graph}%
\author{Dirk Oliver Theis}%
\address{DOT: Fakult\"at f\"ur Mathematik\\
  Otto-von-Guericke-Universit\"at Magdeburg\\
  Universit\"atsplatz~2\\
  39106~Magdeburg\\
  Germany}%
\email{theis@ovgu.de}%
  \date{Wed Aug 31 18:40:28 CEST 2011}
\begin{abstract}
  The biclique covering number of a bipartite graph~$G$ is the minimum number of complete bipartite subgraphs (bicliques) whose union contains every edge of~$G$.
  In this little note we compare three lower bounds on the biclique covering number: A bound $\jk(G)$ proposed by Jukna \& Kulikov (Discrete Math.\ 2009); the
  well-known fooling set bound $\fool(G)$; the ``tensor-power'' fooling set bound $\fool^\infty(G)$.  We show
  \begin{equation*}
    \jk \le \fool \le \fool^\infty \le \min_Q (\rk Q)^2,
  \end{equation*}
  where the minimum is taken over all matrices with a certain zero/nonzero-pattern.  Only the first inequality is really novel, the third one generalizes a result of 
  Dietzfelbinger, Hromkovi{\v{c}}, Schnitger (1994).
  \\%%
  We also give examples for which $\fool \ge (\rk)^{\log_4 6}$ improving on Dietzfelbinger et al.

  \smallskip\noindent%
\end{abstract}
\maketitle

% \begin{itemize}
% \item Check if Lemma~1 really isn't known!\todo{!!}
% \item Check conjecture makes sense!\todo{!!}
% \end{itemize}

\section{Introduction}

For a bipartite graph~$G$, the minimum number of complete bipartite subgraphs (bicliques) whose union contains every edge of~$G$, is called the \textit{biclique
  covering number} $\bc(G)$ of~$G$.  For a 01-matrix~$M$, the minimum number of all-1s submatrices covering all 1s in~$M$ is called the \textit{rectangle covering
  number} of~$M$; the \textit{nondeterministic communication complexity}~\cite{KushilevitzNisan97} of~$M$ is the base-2 logarithm of the rectangle covering number.
If~$M$ is considered as a matrix over the Boolean semiring (i.e., the semiring with two elements $0$ and~$1$, with the usual multiplication and the maximum as
addition), then the minimum number~$q$ for which there exist a matrix $B$ with~$q$ columns and a matrix~$C$ with~$q$ rows such that $M = BC$ (over the Boolean
semiring) is the \textit{Boolean rank}~\cite{CaenGregoryPullman81} of~$M$.

These three concepts all define the same quantity, where one passes between a bipartite graph~$G$ and its bipartite adjacency matrix~$M(G)$.  By abuse of notation,
we will identify between 01-matrices and bipartite graphs.  For example, we will write $\bc(M)$ for a 01-matrix~$M$.

Determining the biclique covering number of a bipartite graph is an NP-hard problem~\cite{Orlin77}.  To make use of it in its many applications, one usually requires
lower bounds.  The probably best-known lower bound is the so-called fooling set bound: if $F$ is a set of edges of~$G$ such that no two of them induce a complete
bipartite graph, then $\bc(G) \ge \sabs{F}$.  (For two edges, to not induce a complete bipartite graph is equivalent to (a)~not being incident and (b)~not being
contained in a cycle of length four.)  The set~$F$ is called a \textit{fooling set.}  The maximum cardinality of a fooling set in~$G$ is denoted by $\fool(G)$ and
known as the fooling set bound.

Recently, Jukna \& Kulikov~\cite{JuknaKulikov09} proved a new lower bound on the biclique covering number.  Their construction is as follows.  Let~$H$ be an induced
subgraph of~$G$ with the property that~$H$ has a perfect matching.  Then
\begin{equation}\label{eq:juku-ieq}
  \bc(G) \ge \tfrac{\sabs{H}^2}{4\Nm{H}},
\end{equation}
where $\sabs{H}$ stands for the number of vertices of~$H$ and $\Nm{H}$ for the number of edges in~$H$.  (Here, and throughout, we follow
Diestel's~\cite{Diestel-GTIII} notation.)  We denote the minimum of the right hand side of~\eqref{eq:juku-ieq} over all~$H$ by $\jk(G)$.
In some situations, Jukna \& Kulikov's bound improves on another well-known lower bound: If $E$ is a matching in~$G$, and $s$ is the largest number of cardinality of
a subset of edges of~$E$ which is contained in a biclique of~$G$, then $\bc(G) \ge \nfrac{\abs{E}}{s}$.  (This is a special case of the so-called ``generalized
fooling sets'' of Dietzfelbinger, Hromkovi{\v{c}}, and Schnitger~\cite{DietzfelbingerHromkovicJurajSchnitger96}.)  Jukna \& Kulikov's bound is better, if~$H$ has few
edges.

A known property of the Boolean rank of a matrix is that $\bc(M^{\otimes j}) \le \bc(M)^j$~\cite{GregoryPullman83,Watts01}, where ``$\otimes j$'' stands for the
$k$th tensor power.  This immediately gives the inequality
\begin{equation*}
  \fool^\infty(M) := \limsup_j \fool( M^{\otimes j})^{\nfrac{1}{j}} \le \bc(M).
\end{equation*}

In this tiny little short communication, we prove the following.

\begin{theorem*}
  Let~$G$ be a bipartite graph with bipartition $U\uplus V$.  Then
  \begin{equation*}
    \jk(G) \le \fool(G) \le \fool^\infty(G) \le \min_{Q} (\rk Q)^2,
  \end{equation*}
  where the minimum extends over all fields~$\kk$ and all $(U\times V)$-matrices $Q$ over~$\kk$ with
  \begin{equation}\label{eq:supp-cond}
    Q_{uv} \ne 0 \quad\text{iff}\quad uv \in E(G).
  \end{equation}
\end{theorem*}

\medskip
\begin{remark}
  If~$H$ is an induced subgraph of~$G$ having a perfect matching, then a fooling set of cardinality at least the right hand side of~\eqref{eq:juku-ieq} can be found
  as a subset of any perfect matching of~$H$.
\end{remark}

\begin{remark}\label{rem:support-graph-vs-mtx}
  In terms of matrices, the condition on the support of~$Q$ reads as follows:
  \begin{equation*}
    Q_{k,\ell} \ne 0 \quad\text{iff}\quad M_{k,\ell} \ne 0.
  \end{equation*}
\end{remark}

% In this short communication, we also propose an example showing that it is tight.
% \begin{proposition}
%   There exists a sequence of bipartite graphs $G_k$, $k=1,2,3,\dots$, for which $\rk M(G_k) = 5^k$ but $\bc(G_k) = (\sqrt 5 - o(1))^k$.
% \end{proposition}

The inequality $\fool(G) \le \bigl( \rk M(G) \bigr)^2$ is due to Dietzfelbinger et al.~\cite{DietzfelbingerHromkovicJurajSchnitger96}, where an infinite family of
graphs/matrices $(G_k)_k$ is given for which $\fool(G_k) \ge \rk(M(G_k))^{\log_3 4}$.  Here we propose a family of examples which improves the exponent.

\begin{proposition}
  There exists an infinite family of bipartite graphs~$(G_k)_k$ for which
  \begin{equation*}
    \fool(G_k) \ge \rk(M(G_k))^{\log_4 6}.
  \end{equation*}
\end{proposition}

The improvement is modest: $\log_3 4 < 1.262$ whereas $\log_4 6 > 1.2924$.

We conjecture the following.

\begin{conjecture*}
  For all bipartite graphs, we have
  \begin{equation*}
    \bc(G) \le \max_Q (\rk Q)^2 \,\log\sabs{G} %%
  \end{equation*}
  where the maximum extends over all real matrices satisfying~\eqref{eq:supp-cond}.
\end{conjecture*}

%%%%%%%%%%%%%%%%%%%%%%%%%%%%%%%%%%%%%%%%%%%%%%%%%%%%%%%%%%%%%%%%%%%%%%%%%%%%%%%%%%%%%%%%%%%%%%%%%%%%%%%%%%%%%%%%%%%%%%%%%%%%%%%%%%%%%%%%%%%%%%%%%%%%%%%%%%%%%%%%
\section{Proof of the theorem}

For a bipartite graph~$G$, define the following graph~$X(G)$.
The vertex set of $X(G)$ equals the edge set of~$G$; two vertices of~$X(G)$ are adjacent, iff they do not induce a complete bipartite subgraph of~$G$.

The following very simple fact certainly must be known, but we have not found it anywhere in the literature.  We give a proof here for the sake of completeness.  The
symbols~$\chi$ and~$\omega$ stand for the chromatic and clique numbers, respectively, of a graph, and $\Theta$ is the Shannon Capacity of a graph.

\begin{lemma}\label{lem:chi}
  For every bipartite graph~$G$,
  \begin{enumerate}[(a)]
  \item\label{lem:chi:chi} $\displaystyle \bc(G) = \chi(X(G))$
  \item\label{lem:chi:omega} $\displaystyle \fool(G) = \omega(X(G))$
  \item\label{lem:chi:Theta} $\displaystyle \fool^\infty(G) = \Theta(\widebar{X(G)})$
  \end{enumerate}
  Moreover, if~$F$ is a fooling set in~$G$, then $F$ is a clique in~$X(G)$.
\end{lemma}
\begin{proof}
  Let $U,V$ be a bipartition for~$G$.
  \begin{enumerate}[(a)]
  \item If $U'\subset U$ and $V'\subset V$ are such that the graph $G[U'\cup V']$ induced by $U'\cup V'$, is a biclique, then, clearly, $U'\times V'$ is an independent
    set in~$X(G)$.  Thus, a covering of the edges of~$G$ by bicliques gives rise to a covering of the vertices of~$X(G)$ by independent sets.
    On the other hand, if~$A$ is an independent set in~$X(G)$, then the subgraph of~$G$ induced by the edge set~$A$ is a biclique.  Hence, if $A_1,\dots,A_k$ is a
    covering of the vertices of~$X(G)$ by independent sets, then there is a biclique covering of~$G$ consisting of at most~$k$ bicliques.
  \item Immediately from the definitions.
  \item Recall the definition of the Shannon Capacity of a graph~$X$ via the strong graph product~$\boxtimes$:
    \begin{equation*}
      \Theta(\bar X) = \lim_j \lt( \omega\lt(  \widebar{ \bigl( \bar X \bigr)^{\boxtimes j} }  \rt) \rt)^{\!\nfrac{1}{j}},
    \end{equation*}
    where $\bar\cdot$ denotes the complement of a graph, and $\boxtimes j$ refers to the $j$-fold strong product: The graph $Y^{\boxtimes j}$ has vertex set
    $V(Y)^j$, and $u,v\in V(Y)^j$ are adjacent if, for all~$i$, $u_i=v_i$ or $u_i\sim v_i$.  The claim~(\ref{lem:chi:Theta}) now follows from~(\ref{lem:chi:omega})
    and the fact
    \begin{equation*}
        X(M^{\otimes j})  = \widebar{   \Bigl( \widebar{X(M)} \Bigr)^{\boxtimes j}   },
    \end{equation*}
    whose verification is straight forward.
  \end{enumerate}
\end{proof}

\paragraph{The following facts are worth noting, too.}
\begin{enumerate}[{$\bullet$}]
\item The generalized fooling set bound of Dietzfelbinger et al.~\cite{DietzfelbingerHromkovicJurajSchnitger96} is the same as the hereditary independence ratio
  bound on the chromatic number
  \begin{equation*}
    \bc(G) \ge \imath(X(G)) := \min \Bigl\{ \frac{\sabs U}{\alpha(U)} \Bigm| U\subset V(X(G)) \Bigr\},
  \end{equation*}
  where $\alpha(U)$ stands for the independence number of the subgraph of~$X(G)$ induced by the vertex set~$U$.
\item The fractional covering number $\bc^*(G)$ of~$G$ (e.g.,~\cite{KarchmerKushilevitzNisan95}) corresponds to the fractional chromatic number of~$X(G)$.
  Surprisingly, this observation yields the following slight improvement over Corollary~2.3 in \cite{KarchmerKushilevitzNisan95}: By standard facts on the
  fractional chromatic number (e.g., Theorem~64.13 in \cite{SchrijverBk03}), we have $\bc(G) \le (1+\ln s)\bc^*(G)$, with~$s$ the number of edges in the largest
  biclique in~$G$.  The improvement is from a factor of $O(\sabs{G})$ to a factor of $1+\ln s$, which is useful if $s \ll \sabs{G}$.
\item Applying known tools from graph coloring gives the following ``new'' (cf.~\cite{Lovasz79}) lower bound on the biclique covering number.  With $n := \abs{G}$,
  let $d \in \RR^n\setminus(0)$ and for every edge $e$ of~$G$, let $x_e \in \RR^n\setminus(0)$.  If two edges $e$ and~$f$ do not induce a biclique, then we require
  that $x_e$ and~$x_f$ are orthogonal.  Then we have
  \begin{equation}\label{eq:theta-bd}
    \bc(G) \ge \sum_{e\in E(G)} \tfrac{1}{\Nm{d}\Nm{x_e}} d^\Tp x_e,
  \end{equation}
  where $\Nm{\cdot}$ stands for the usual Euclidean norm.  The maximum of the right hand side of~\eqref{eq:theta-bd} over all choices of~$d$ and $(x_e)_e$ with said
  properties, can be computed in polynomial time, and is always at least as large as the fooling set bound.
\end{enumerate}

\paragraph{For the proof of the theorem,}  %%
we also need the following fact, which is a \textit{very} slight generalization of Lemma~2.9 in~\cite{DietzfelbingerHromkovicJurajSchnitger96}.

\begin{lemma}[Dietzfelbinger et al.~\cite{DietzfelbingerHromkovicJurajSchnitger96}]\label{lem:fool-rk-ieq}%
  Let~$Q$ be as in the theorem.  Then
  \begin{equation*}
    \fool(G) \le (\rk Q)^2.
  \end{equation*}
\end{lemma}
The proof in~\cite{DietzfelbingerHromkovicJurajSchnitger96} can be copied almost word for word.  For the sake of completeness, we give a here a self-contained
version.
\begin{proof}
  Let~$F$ be a fooling set in~$G$, and~$Q$ as described.  W.l.o.g., we may assume that~$F$ spans~$G$.  (Otherwise take the corresponding subgraph of~$G$, which
  corresponds to taking a sub-matrix of~$Q$.)  We may also assume that $U=V=\{1,\dots,n\}$, so that~$Q$ is an $(n\times n)$-matrix.  Let $Q = XY$ be a rank
  factorization of~$Q$, i.e., with $r := \rk Q$, $X$ is an $(n\times r)$-matrix and $Y$ is an $(r\times n)$-matrix.  Let $x_1,\dots,x_n \in\kk^r$ be the rows of~$X$,
  and $y_1,\dots,y_n\in\kk^r$ the columns of~$Y$, so that $Q_{k,\ell} = \ip{x_k,y_\ell}$, where~$\ip{\cdot,\cdot}$ denotes the inner product.  The tensor product
  $\kk^r \otimes \kk^r$ has an inner product satisfying $\ip{ \xi_1\otimes \xi_2 , \eta_1\otimes \eta_2 } = \ip{\xi_1,\eta_2}\ip{\xi_2,\eta_1}$.  In this inner
  product space of dimension~$r^2$, the tensors $x_k\otimes y_k$, $k=1,\dots,n$, form an orthogonal system (i.e., the inner product of an element with itself is
  non-zero, but the inner product of any two distinct elements is zero), implying $n \le r^2$.
\end{proof}

\paragraph{Now we complete the proof of the theorem.}
\begin{enumerate}[(1)]
\item $\jk(G) \le \fool(G)$.\\
  Let~$H$ be as in the theorem, and $r := \sabs{H}/2$.  Choose an arbitrary perfect matching of~$G$.  We may assume that, for every edge $uv$ of the perfect
  matching, we have for the number of neighbors $d_{uv} := \deg(u) + \deg(v) -2$ of $u$ and $v$,
  \begin{equation*}
    \frac{ r^2 }{ \Nm{H} } \le \frac{ 2r-1 }{ d_{uv}+1 } < \frac{ 2r }{ d_{uv} },
  \end{equation*}
%   \begin{equation*}
%     \frac{  (r-1)^2 }{\Nm{H} - (d_{uv}+1)}
%     = 
%     \frac{  r^2 - (2r-1) }{\Nm{H} - (d_{uv} +1)}
%     >
%     \frac{ r^2 }{ \Nm{H} }
%   \end{equation*}
%   if 
%   \begin{equation*}
%     \frac{ r^2 }{ \Nm{H} } > \frac{ 2r-1 }{ d_{uv}+1 }.
%   \end{equation*}
  because, if the left inequality did not hold, then deleting $u$ and~$v$ and their incident edges would improve the bound.  In particular, we have $d_{uv} \le
  \nfrac{2\Nm H}{r}$.
  
  Consider now the subgraph $X'$ of~$X(H)$ induced by a the edges of a perfect matching in~$H$.  The graph $X'$ has~$r$ vertices, and every vertex of~$X'$ has at
  least $\nfrac{d_{uv}}{2} = r-\frac{\Nm H}{r}$ neighbors (in~$X'$), so that $\Nm{X'} \ge \frac12(r^2 - \Nm H)$.  By Tur\'an's theorem, $X'$ contains a clique of
  size at least
  \begin{equation*}
    \frac{r^2}{r^2 - 2\Nm{X'}} \ge \frac{r^2}{\Nm H}.
  \end{equation*}
  The statement of the theorem now follows from Lemma~\ref{lem:chi}(\ref{lem:chi:omega}).
  
\item $\fool(G) \le \fool^\infty(G)$.  Obvious from $\omega(\cdot) \le \Theta(\bar \cdot)$ and Lemma~\ref{lem:chi}(\ref{lem:chi:omega},\ref{lem:chi:Theta}).
  
\item $\fool^\infty(G) \le \min_Q (\rk Q)^2$.\\
  Let~$Q$ be a matrix as described.
  Then~$Q^{\otimes j}$ has the property that for $u\in U^j$ and $v\in V^j$, we have
  \begin{eqnarray*}
    ~          & \bigl( Q^{\otimes j} \bigr)_{uv}  = \prod_{i=1}^j Q_{u_j v_j} \ne 0& \\
    \text{iff} & \bigl( M^{\otimes j} \bigr)_{uv} = \prod_{i=1}^j M_{u_j v_j} = \ne 0.
  \end{eqnarray*}
  Put differently, $Q^{\otimes j}$ satisfies the conditions for Lemma~\ref{lem:fool-rk-ieq}.  Hence, since $\rk(Q^{\otimes j}) = (\rk Q)^j$, we conclude that
  \begin{equation*}
    \fool(M^{\otimes j})
    \le
    \bigl(  \rk( Q^{\otimes j} )  \bigr)^{2}
    =
    ( \rk Q )^{2j}
  \end{equation*}
  which implies $\fool^\infty(G) \le \min_Q (\rk Q)^2$, as claimed.
  \qed
\end{enumerate}

\section{Proof of the proposition}

Consider the following matrix.
\begin{equation*}
  M :=
  \begin{pmatrix}
    1 & 1 & 1 & 0 & 0 & 0 \\
    0 & 1 & 1 & 1 & 0 & 0 \\
    0 & 0 & 1 & 1 & 1 & 0 \\
    0 & 0 & 0 & 1 & 1 & 1 \\
    1 & 0 & 0 & 0 & 1 & 1 \\
    1 & 1 & 0 & 0 & 0 & 1 \\
  \end{pmatrix}.
\end{equation*}
Obviously, we have $\fool(M) = 6$.  Since $x_1 := (2,-1,-1,2,-1,-1)^\Tp$ and $x_2 := (-1,2,-1,-1,2,-1)^\Tp$ are two linearly independent vectors with $Mx_i = 0$,
$i=1,2$, we have that $\rk M \le 4$.  Using the arguments from the previous section, we obtain
\begin{equation*}
  \fool M^{\otimes k} \ge 6^k \qquad\text{and}\qquad\rk M^{\otimes k} \le 4^k,
\end{equation*}
and conclude that $\fool M^{\otimes k} \ge \Bigl( \rk M^{\otimes k} \Bigr)^{\log_4 6}$.  The proposition follows when we let~$G_k$ be the bipartite graph whose
bipartite adjacency matrix equals $M^{\otimes k}$.  (We note that the tensor power construction is the same as in~\cite{DietzfelbingerHromkovicJurajSchnitger96};
only the initial matrix~$M$ differs.)

% %%%%%%%%%%%%%%%%%%%%%%%%%%%%%%%%%%%%%%%%%%%%%%%%%%%%%%%%%%%%%%%%%%%%%%%%%%%%%%%%%%%%%%%%%%%%%%%%%%%%%%%%%%%%%%%%%%%%%%%%%%%%%%%%%%%%%%%%%%%%%%%%%%%%%%%%%%%%%%%%
% \section{Acknowledgments}
%
% The author would like to thank the \textit{Casa del Gelato} in Magdeburg for making quite good Espresso.

\providecommand{\bysame}{\leavevmode\hbox to3em{\hrulefill}\thinspace}
\providecommand{\MR}{\relax\ifhmode\unskip\space\fi MR }
% \MRhref is called by the amsart/book/proc definition of \MR.
\providecommand{\MRhref}[2]{%
  \href{http://www.ams.org/mathscinet-getitem?mr=#1}{#2}
}
\providecommand{\href}[2]{#2}

\end{document}